\documentclass[a4paper,11pt,onecolumn]{amsart}
\usepackage{amsfonts}
\usepackage{pb-diagram}
\title[Actions on orderings of groups]{Faithful actions of automorphisms on the space of orderings of a group}

\newtheorem{thm}{Theorem}
\newtheorem{lemma}[thm]{Lemma}

\newtheorem{prop}[thm]{Proposition}

\numberwithin{thm}{section}

\newcommand{\al}{\alpha}

\newcommand{\gam}{\Gamma}

\newcommand{\bZ}{\mathbb{Z}}
\newcommand{\bR}{\mathbb{R}}

\newcommand{\bN}{\mathbb{N}}

\DeclareMathOperator{\Aut}{Aut}
\DeclareMathOperator{\Out}{Out}

\DeclareMathOperator{\Homeo}{Homeo}
\DeclareMathOperator{\Mod}{Mod}

\newcommand{\bP}{\mathbb{P}}

\newcommand{\bQ}{\mathbb{Q}}

\newcommand{\mP}{\mathcal{P}}

\author[T. Koberda]{Thomas Koberda}
\address{Department of Mathematics\\ Harvard University\\ 1 Oxford St.\\ Cambridge, MA 02138 }
\email{ koberda@math.harvard.edu}
\begin{document}
\begin{abstract}
In this article we study the space of left-- and bi--invariant orderings on a torsion--free nilpotent group $G$.  We will show that generally the set of such orderings is equipped with a faithful action of the automorphism  group of $G$.  We prove a result which allows us to establish the same conclusion when $G$ is assumed to be merely residually torsion--free nilpotent.  In particular, we obtain faithful actions of mapping class groups of surfaces.  We will draw connections between the structure of orderings on residually torsion--free nilpotent, hyperbolic groups and their Gromov boundaries, and we show that in those cases a faithful $\Aut(G)$--action on the boundary is equivalent to a faithful $\Aut(G)$ action on the space of left--invariant orderings.
\end{abstract}
\maketitle
\begin{center}
\today
\end{center}
\tableofcontents
\section{Introduction}
The purpose of this article is to show that the space of left--invariant orderings of a residually torsion--free nilpotent group $G$ is sufficiently rich as to admit a faithful action of $\Aut(G)$.

Let $G$ be a finitely generated group.  A fundamental and often quite difficult problem in the combinatorial group theory of $G$ is to describe the space of orderings on $G$.  A {\bf left--invariant ordering} on $G$ is a relation $\leq$ on $G$ which is a total ordering on the elements of $G$, together with the following left--invariance property: for all triples $a,b,c\in G$, $a\leq b$ implies $ca\leq cb$.  An ordering is called {\bf right--invariant} if the analogous right--invariance property holds.  An ordering is called {\bf bi--invariant} if it is both left-- and right--invariant.  It is easy to check that an ordering is bi--invariant if and only if it is left--invariant and conjugation--invariant.

Many groups admit no left--invariant orderings at all.  For instance, the presence of torsion precludes orderability.  Some groups admit finitely many orderings, and the book of Botto Mura and Rhemtulla \cite{MR} describes some aspects of the theory of groups with finitely many orderings.

It is sometimes useful to observe that orderings on a group naturally occur in pairs.  For each ordering $\leq$, there is a natural ordering $\leq^{op}$ called the {\bf opposite ordering}, given by $g\leq^{op}h$ if and only if $h\leq g$.  Many groups admit uncountably many orderings.

To organize the set of all orderings of a group, one defines the {\bf space of left--invariant or bi--invariant orderings} on the group, denoted $LO(G)$ in the case of left--invariant orderings and $O(G)$ in the case of bi--invariant orderings.  To define this space and equip it with a good topology, we first define the notion of a {\bf positive cone} $\mP$ of an ordering.  Given an ordering $\leq\in LO(G)$ or $O(G)$, we set \[\mP=\mP(\leq)=\{g\in G \textrm{ such that } 1<g\}.\]  This gives us a canonical bijective correspondence between orderings and  certain subsets of $G$.  Indeed, to recover an ordering, we declare $g<h$ if and only if $g^{-1}h\in\mP$.

In order for a subset of $G$ to be the positive cone of some left--invariant ordering, it must satisfy some axioms:
\begin{enumerate}
\item
$\mP\cup\mP^{-1}=G\setminus \{1\}$, where $\mP^{-1}$ denotes the set of inverses of elements of $\mP$.
\item
$\mP\cap\mP^{-1}=\emptyset$.
\item
$\mP\cdot\mP\subset\mP$.
\end{enumerate}

$\mP$ will be the positive cone of some bi--invariant ordering if in addition $\mP$ is $G$--conjugation invariant.

The power set of subsets of $G$ comes with a natural topology which gives it the structure of a Cantor set.  Precisely, two subsets are close if they agree on large finite subsets.  This Cantor set will be metrizable whenever $G$ is countable.  In particular, the power set of $G$ can be viewed as \[\{0,1\}^G,\] where the two point set has the discrete topology and the product has the product topology.  Two points in the power set of $G$ are close in this topology if they agree on a large finite subset.

It is not difficult to show that the conditions for a set to be the positive cone of a left-- or bi--invariant ordering are closed conditions in the natural topology on the power set.  The details of the proof can be found in Chapter $14$ of the book \cite{DDRW} by Dehornoy, Dynnikov, Rolfsen and Wiest.  Thus, $LO(G)$ and $O(G)$ can be viewed as closed subsets of a Cantor set.  The topology of this space for various groups has been studied by various authors, such as by Navas for free groups in in \cite{N}, by Navas and Rivas for Thompson's group $F$ in \cite{NR}, and by Sikora for finitely generated torsion--free abelian groups in \cite{S}.

The groups $\Aut(G)$ and $\Out(G)$ both have natural actions on $LO(G)$ and $O(G)$ respectively.  $G$ acts on $LO(G)$ by conjugation, so $\Out(G)$ also acts on the $G$--conjugation orbits in this space.  These actions are given by pulling back and ordering $\leq$ to an ordering $\leq_{\phi}$ via the automorphism $\phi$.  Precisely, we define $g\leq_{\phi} h$ if and only if $\phi(g)\leq\phi(h)$.  In the case of the action of $\Out(G)$, the conjugation action of $G$ on $O(G)$ is trivial, so it does not matter which automorphism representative for an outer automorphism we choose.  It is easy to check that the two actions are by homeomorphisms.  One sees that this way we get maps \[\psi_a:\Aut(G)\to\Homeo(LO(G))\] and \[\psi_o:\Out(G)\to\Homeo(O(G)).\]  

These maps, particularly the first, are the primary focus of this paper.  Recall that a group $G$ is called {\bf residually torsion--free nilpotent} if every non--identity element of $G$ persists in some torsion--free nilpotent quotient of $G$.  Examples of residually torsion--free nilpotent groups include free groups, surface groups, right-angled Artin groups and pure braid groups.  With this terminology, we can state the main result of this paper:

\begin{thm}\label{t:main}
Let $G$ be a finitely generated, residually torsion--free nilpotent group.  Then the map \[\psi_a:\Aut(G)\to\Homeo(LO(G))\] is injective.
\end{thm}

In particular, the conclusions of Theorem \ref{t:main} hold for mapping class groups of surfaces (with a marked point) and automorphism groups of free groups.  Theorem \ref{t:main} shows that there are many essentially different positive cones in residually torsion--free nilpotent groups which are not preserved by automorphisms of the group.

The proof of Theorem \ref{t:main} is of a very similar flavor to the proof of asymptotic linearity of the mapping class group, one of the principal results in \cite{K}.  Asymptotically faithful actions of mapping class groups have been of recent interest to various authors, such as Andersen in \cite{A}.

As an amplification of the ideas of Theorem \ref{t:main}, we will show that when $G$ is residually torsion--free nilpotent and hyperbolic, $LO(G)$ recovers the boundary $\partial G$.  We will be able to show:
\begin{thm}\label{t:boundary}
Suppose that $G$ is residually torsion--free nilpotent and hyperbolic.  Then $\Aut(G)$ acts faithfully on $\partial G$.
\end{thm}

In the case that $G$ is a surface group, Theorem \ref{t:boundary} can be viewed as a generalization of the classical result of Nielsen, namely that the mapping class group $\Mod_{g,1}$ of a surface of genus $g\geq 2$ with one marked point acts faithfully on the circle.  For more details, consult the book of Casson and Bleiler \cite{CB}.  It seems that there were few if any connections between orderings on groups and geometric group theory appearing anywhere in literature.  It thus appears that Theorem \ref{t:boundary} gives an example of such a connection.

It is unlikely that one can easily remove the residual condition on $G$ in the statement of Theorem \ref{t:boundary}, since hyperbolic groups can be so diverse.  It is not even known whether or not every hyperbolic group is residually finite or virtually torsion--free.  For some discussion of virtual properties of hyperbolic groups, the reader might consult the paper \cite{KW} of I. Kapovich and D. Wise.

We close the paper by showing that Theorem \ref{t:main} does not hold in general:

\begin{prop}\label{p:klein}
Let $K$ be the fundamental group of the Klein bottle.  Then $\Aut(K)$ does not act faithfully on $LO(K)$ and $\Out(K)$ does not act faithfully on $O(K)$ nor on conjugacy classes of elements of $LO(K)$.
\end{prop}

\section{Acknowledgements}
This paper benefitted from conversations with M. Bestvina, B. Farb, T. Church, C. McMullen, C. Taubes and B. Wiest.  The author thanks P. Hubert for asking whether orderings and Gromov hyperbolicity are related.  The author finally thanks the referees for careful reading and useful comments and corrections, and for contributing some simplification to the proofs.  The author was supported by an NSF Graduate Research Fellowship for part of the time that this research was carried out.

\section{Abelian groups}
In order to prove Theorem \ref{t:main}, we will need to understand the conclusion of the theorem for finitely generated torsion--free abelian groups.  Our goal is to prove:

\begin{lemma}\label{t:abelian}
$GL_n(\bZ)$ acts faithfully on $O(\bZ^n)$ under the homomorphism $\psi_a$.
\end{lemma}

First, we must understand the structure of $LO(\bZ^n)=O(\bZ^n)$.  When $n=1$, it is evident that this set has exactly two points.  When $n>1$, Sikora proved in \cite{S} that $O(\bZ^n)$ is a Cantor set.  To adapt Sikora's Theorem to our setup, we will be quite explicit about a construction of certain orderings on $\bZ^n$.

We begin by identifying some useful orderings on $\bZ^n$.  Let $Z$ denote an rational hyperplane in $\bR^n$.  Then $Z$ will help determine many positive cones on $\bZ^n$ as follows: choose a half of $\bR^n\setminus Z$ to be positive.  Then choose a hyperplane within $Z$ and declare a half of $Z$ to be positive.  Continuing this process, we eventually declare each nonzero integral point in $\bR^n$ to be either positive or negative.  It is easy to see that we in fact obtain a positive cone on $\bZ^n$ this way.

It follows that a flag of rational subspaces of $\bR^n$ together with a choice of half--space in each dimension gives rise to an ordering on $\bZ^n$.  We will call orderings which arise in this fashion {\bf flag orderings}.

Note that if $Z$ is an irrational hyperplane in the sense that it contains no rational points other than the origin, $Z$ automatically already determines exactly two orderings: one for each choice of positive halfspace.

We have two perspectives on orderings of $\bZ^n$.  One comes from choosing irrational hyperplanes and rational flags as above, and the other comes the definition of a positive cone.  It is not immediately clear how to reconcile these two descriptions of the orderings on $\bZ^n$, even in the case $n=2$.  When $n=2$, we have a map from $O(\bZ^2)$ to $\bR\bP^1$.  This map is given by sending an ordering to the line which separates the positive half--plane from the negative half.  The fiber over an irrational point in $\bR\bP^1$ consists of two points, one for each choice of positive half--plane.  The fiber over a rational point consists of four points, corresponding to the two choices for positive half--plane and the two choices for positive half--line.  Thus, one can see that the space of orderings should not be considered with an analytic topology, but rather with a topology which more closely resembles a totally disconnected one.

The exact nature of the topology on $LO(G)$ and $O(G)$ is not important for the purposes of this article and we will not discuss the topology much further, other than to remark that the ``fibration" $O(\bZ^2)\to\bR\bP^1$ is continuous in the Cantor set topology on $O(\bZ^2)$ and the usual topology on $\bR\bP^1$.  A discussion of this map can be found in \cite{S}.

The rational flag orderings occupy a special place in the study of orderings on $\bZ^n$, since they are dense in the space of all orderings on $\bZ^n$:

\begin{lemma}\label{t:halfspace}
Let $V=\{v_1,\ldots,v_{n+1}\}\subset\bZ^n$ be nonzero vectors which do not lie within a single closed rational half space in $\bR^n$.  Let $S$ denote the semigroup generated by these vectors.  Then $0\in S$.
\end{lemma}
\begin{proof}
The conclusion for $n=1$ is trivial.
For the general case, fix a basis $\{x_1,\ldots,x_{n+1}\}$ for $\bQ^{n+1}$ and let \[A:\bQ^{n+1}\to\bQ^n\] be the linear map which sends the vector $x_i$ to the vector $v_i$.  The map $A$ evidently has a nontrivial kernel.

Suppose that the conclusion of the lemma fails.  Then there is an integral vector \[w=(a_1,\ldots,a_{n+1})\] contained in the kernel of $A$ such that neither $w$ nor $-w$ is contained in the closed positive orthant of $\bQ^{n+1}$.  But then there are two indices, which we may assume are $1$ and $2$, with $a_1>0$ and $a_2<0$.  Then we have \[\sum_{i=3}^{n+1}a_iv_i+a_1v_1=-a_2v_2,\] so that $v_1$ and $v_2$ are on the same side of the hyperplane spanned by $\{v_3,\ldots,v_{n+1}\}$.  But then the vectors $\{v_1,\ldots,v_{n+1}\}$ are all contained in one closed halfspace, a contradiction of the hypotheses.
\end{proof}

The author is indebted to the referee for simplifying the proof of Lemma \ref{t:halfspace}.

\begin{prop}\label{p:dense}
The set of flag orderings on $\bZ^n$ is dense in the space of orderings on $\bZ^n$ in the Cantor set topology.
\end{prop}
\begin{proof}
Let $\{(a_1,b_1),\ldots,(a_m,b_m)\}$ be a collection of pairs of distinct lattice points and let $P\in O(\bZ^n)$.  Suppose that according to $P$, $a_i<b_i$ for all $i$.  We will show that there is a flag ordering in which these relations also hold.  This will imply that in any open subset of $O(\bZ^n)$ containing $P$, there is a flag ordering.

By definition, $(b_i-a_i)\in P$ for each $i$.  By Lemma \ref{t:halfspace}, all the vectors $\{(b_i-a_i)\}$ must lie in a closed rational halfspace.  If there is a rational hyperplane $H$ such that all the vectors $\{(b_i-a_i)\}$ are in one open half space defined by $H$, then we are done.  Otherwise, we consider the elements of $\{(b_i-a_i)\}$ which lie in $H$.  A repeated application of Lemma \ref{t:halfspace} shows that there is a flag ordering on $\bZ^n$ where all the vectors $\{(b_i-a_i)\}$ are positive.
\end{proof}

We are now ready to prove Lemma \ref{t:abelian}.

\begin{proof}[Proof of Lemma \ref{t:abelian}]
We claim that in fact $GL_n(\bZ)$ acts faithfully on the set of flag orderings of $\bZ^n$.  Let $1\neq A\in GL_n(\bZ)$ be an automorphism which preserves every flag ordering on $\bZ^n$.  Then $A$ must preserve each rational hyperplane in $\bR^n$.  Indeed, if $H$ and $J$ are distinct rational hyperplanes then $H$ and $J$ cut $\bR^n$ into halfspaces \[\{S_{1,H},S_{2,H},S_{1,J},S_{2,J}\}\] whose intersections with $\bZ^n$ are all different.  Therefore if $A$ sends $H$ to $J$ then $A$ acts nontrivially on $O(\bZ^n)$.

It follows that $A$ preserves each rational hyperplane and therefore acts trivially on $\bP^{n-1}(\bQ)$ (via the dual action).  It follows that $A$ is trivial in $\bP GL_n(\bZ)$ and is therefore a scalar multiple of the identity.  If $A$ is nontrivial and integral then it would have to be $-I$.  It is clear that $-I$ acts nontrivially on $O(\bZ^n)$.
\end{proof}

\section{Extensions and pullbacks of orderings}
Other than the machinery of orderings on abelian groups, certain extension and pullback theorems for orderings on torsion--free nilpotent groups will be very important for the proof of Theorem \ref{t:main}.  Up to this point in our discussion of orderings on groups, we have been considering positive cones which contain ``half" of the nonidentity elements in a group.  If we are given a positive cone $\mP$ which is partial in the sense that $\mP\cup\mP^{-1}$ is properly contained in $G\setminus\{1\}$, we call $\mP$ a {\bf partial ordering}.  A partial ordering $\mP$ is bi--invariant if it is conjugation--invariant.  We now quote the following two strong theorems, the first due to Rhemtulla in \cite{Rh} and the second due to Mal'cev in \cite{M} (see also \cite{MR}):

\begin{thm}
Let $N$ be a finitely generated torsion--free nilpotent group and $P$ a partial ordering on $N$.  Then $P$ extends to a total ordering on $N$.
\end{thm}

\begin{thm}
Let $N$ be a finitely generated torsion--free nilpotent group and $P$ a bi--invariant partial ordering on $N$.  Then $P$ extends to a total bi--invariant ordering on $N$.
\end{thm}

Mal'cev actually proved that it suffices for $N$ to be locally torsion--free nilpotent.

The two extension theorems above can be restated as follows:

\begin{thm}
Let $N$ be a torsion--free nilpotent group and let $\{1\}\neq N'<N$ be a subgroup.  Then the restriction map \[\rho_L:LO(N)\to LO(N')\] is surjective.  If in addition $N'$ is normal then the restriction map \[\rho_B:O(N)\to O(N')^N\] is surjective, where $O(N')^N$ denotes the $N$--invariant bi--invariant orderings on $N'$.
\end{thm}

We will call the previous results {\bf Rhemtulla's} and {Mal'cev's Extension Theorems}, respectively.

We will often encounter a situation where $N$ is a torsion--free nilpotent quotient of a group $G$ equipped with an ordering, and we wish to produce an ordering of $G$ which is compatible with the quotient map $G\to N$ and the given ordering on $N$.

Let $G$ be a finitely generated group.  We will write $\{\gamma_i(G)\}$ for the lower central series of $G$.  This series is defined by $\gamma_1(G)=G$ and $\gamma_{i+1}(G)=[G,\gamma_i(G)]$.  A group is called {\bf residually nilpotent} if \[\bigcap_{i>0}\gamma_i(G)=\{1\}.\]  The usual definition of residual nilpotence says that $G$ is residually nilpotent if for each nonidentity $g\in G$, there exists a nilpotent quotient $N_g$ of $G$ where $g$ is not mapped to the identity.  These two definitions are equivalent.  Indeed, if $N_g$ satisfies $\gamma_i(N_g)=\{1\}$ and is a quotient of $G$, then $N_g$ is a quotient of the group obtained from $G$ by declaring $\gamma_i(G)=1$.  Conversely, if $g\neq 1$ then there is some $i$ for which $g\notin\gamma_i(G)$, whence $g$ survives in $G/\gamma_i(G)$.

A finitely generated group is called {\bf residually torsion--free nilpotent} if for each nontrivial $g\in G$, there is an $i$ for which $g$ maps to an infinite order element of $G/\gamma_i(G)$.  Again, the usual definition of a residually torsion--free nilpotent group $G$ says that each nonidentity $g\in G$ survives in a torsion--free nilpotent quotient $T_g$ of $G$.  If $T_g$ is a quotient of $G/\gamma_i(G)$ (which it must for some $i$, since $T_g$ is nilpotent) then the image of $g$ in $G/\gamma_i(G)$ has infinite order.  Conversely, the existence of a torsion--free quotient of $G/\gamma_i(G)$ follows from the following well--known fact about nilpotent groups:

\begin{lemma}
Let $N$ be a finitely generated nilpotent group.
\begin{enumerate}
\item
The elements of finite order in $N$ generate a finite normal subgroup $T(N)$.
\item
The quotient $N/T(N)$ is torsion--free.
\end{enumerate}
\end{lemma}
\begin{proof}
See \cite{R}, for instance.
\end{proof}

The previous lemma allows us to modify the lower central series of a residually torsion--free nilpotent group $G$ in a way which will be useful in further discussion.  We will let $\gamma_i^T(G)$ be the kernel of the composition map \[G\to G/\gamma_i(G)\to (G/\gamma_i(G))/T(G/\gamma_i(G)).\]  Then $G/\gamma_i^T(G)$ is torsion--free, and \[\bigcap_{i>0}\gamma_i^T(G)=\{1\}.\]  Observe that $\gamma_i^T(G)$ is characteristic in $G$ and that if $i<j$ then $\gamma_j^T(G)<\gamma_i^T(G)$.

Observe that since for any nilpotent group $N$, the subgroup $T(N)$ is finite, we immediately see that $\gamma_i(G)<\gamma_i^T(G)$ as a finite index subgroup.  It follows that the groups \[\gamma_i(G)/\gamma_{i+1}(G)\] and \[\gamma_i^T(G)/\gamma_{i+1}^T(G)\] are commensurable.

It is a classical result that if $\phi\in\Aut(G)$ acts trivially on the abelianization $G^{ab}$ of $G$ then it also acts trivially on $\gamma_i(G)/\gamma_{i+1}(G)$.  A detailed proof and discussion can be found in \cite{BL}, for instance.  One perspective on this fact is that there is a natural surjective map from the $i$--fold tensor product of $G^{ab}$ to $\gamma_i(G)/\gamma_{i+1}(G)$, given by the commutator bracket.

On the one hand, an automorphism $\phi$ may act nontrivially on $G^{ab}$ and yet descend to the identity on $G/\gamma_2^T(G)$.  On the other hand, we have the following:

\begin{lemma}\label{l:trivial}
Suppose $\phi\in\Aut(G)$ acts trivially on $G/\gamma_2^T(G)$.  Then $\phi$ acts trivially on $\gamma_i^T(G)/\gamma_{i+1}^T(G)$ for all $i$.
\end{lemma}
\begin{proof}
Let $\omega$ be an $i$--fold tensor of elements of $G^{ab}$, viewed as an element of $\gamma_i(G)/\gamma_{i+1}(G)$.  If any factor of the tensor $\omega$ has finite order then multi--linearity of the tensor product implies that the image of $\omega$ in $\gamma_i(G)/\gamma_{i+1}(G)$ has finite order as well.  Thus the natural surjective map \[\bigotimes G^{ab}\to\gamma_i(G)/\gamma_{i+1}(G)\] descends to a natural map \[\bigotimes G/\gamma_2^T(G)\to\gamma_i^T(G)/\gamma_{i+1}^T(G)\] whose image has finite index in the target.  If $\phi$ acts trivially on $G/\gamma_2^T(G)$ then it must act trivially on a finite index subgroup of $\gamma_i^T(G)/\gamma_{i+1}^T(G)$.  Since the latter is torsion--free, it follows that $\phi$ induces the trivial automorphism of $\gamma_i^T(G)/\gamma_{i+1}^T(G)$.  The conclusion follows.
\end{proof}

The following proposition is well--known (a proof with applications to the theory of braid orderings can be found in \cite{DDRW}) but we recall a proof for the convenience of the reader and because the proof will motivate further discussion:

\begin{prop}\label{p:standard}
Let $G$ be finitely generated and residually torsion--free nilpotent.  Then $O(G)$ is nonempty.
\end{prop}
\begin{proof}
For each $i$, write $N_i$ for $G/\gamma_i^T(G)$, and let $Z_i$ denote the kernel of the map $N_i\to N_{i-1}$, where by convention $N_0=\{1\}$.  Note that \[Z_i=\gamma_i^T(G)/\gamma_{i+1}^T(G).\]  Observe that each $Z_i$ is a finitely generated free abelian group, and the conjugation action of $N_i$ on $Z_i$ is trivial.  The reason for the second claim is that $Z_i$ is virtually central in $N_i$, so there can be no action of $N_i$ on $Z_i$ which is nontrivial and yet restricts to the identity on a finite index subgroup.

Choose an arbitrary ordering on each $Z_i$.  We obtain an element of $O(G)$ as follows: let $g,h\in G$.  Suppose $N_i$ is the first such quotient of $G$ in which $g^{-1}h$ survives.  Then by minimality of $i$, we have $g^{-1}h\in Z_i$ under the map $G\to N_i$.  If the image of $g^{-1}h$ is positive in $Z_i$, we declare $g<h$ in $G$.  It is easy to see that this defines a bi--invariant ordering on $G$.
\end{proof}

Orderings as in Proposition \ref{p:standard} are called {\bf standard orderings}.  Using ideas similar to those in the proof of Proposition \ref{p:standard}, we can pull back orderings on torsion--free nilpotent quotients of a residually torsion--free nilpotent group in a way which we call the {\bf standard ordering construction}.

\begin{lemma}\label{l:soc}
Let $G$ be a residually torsion--free nilpotent group and let $N=G/\gamma_i^T(G)$.  Suppose we are given an ordering $\overline{P}\in LO(N)$.  Then there exists an ordering $P\in LO(G)$ which is a pullback of $\overline{P}$ in the following sense:
for all $\overline{g},\overline{h}\in N$ and any preimages $g,h\in G$, we have $g<h$ in $G$ if and only if $\overline{g}<\overline{h}$ in $N$.
\end{lemma}
\begin{proof}
For each $j>i$, choose an arbitrary ordering on $Z_j$.  Let $g,h\in G$.  If $g^{-1}h$ is nontrivial in $N$ then we declare $g<h$ if and only if $g^{-1}h\in P$ under the projection $G\to N$.  Otherwise we may find, as in the proof of Proposition \ref{p:standard}, a minimal $j$ for which $g^{-1}h$ survives in some $Z_j$.  We declare $g<h$ if the image $g^{-1}h$ is positive in the ordering on $Z_j$.
\end{proof}

\section{Representations of automorphism groups and the boundary of a hyperbolic group}
In this section, we prove Theorem \ref{t:main}.

\begin{thm}\label{t:m1}
Let $G$ be a residually torsion--free nilpotent group and let $1\neq \phi\in\Aut(G)$.  Then $\phi$ acts nontrivially on $LO(G)$.
\end{thm}
\begin{proof}
Clearly we may suppose that $\phi$ acts trivially on $G/\gamma_2^T(G)$, since otherwise we may choose an ordering on $G/\gamma_2^T(G)$ which is not preserved by $\phi$ by Lemma \ref{t:abelian}, and then pull it to all of $G$ by a standard ordering construction, as in Lemma \ref{l:soc}.

Suppose $\phi$ acts trivially on $G/\gamma_2^T(G)$ but that $\phi$ is a nontrivial automorphism of $G$.  Let $i$ be minimal so that $\phi$ acts nontrivially on $N_i=G/\gamma_i^T(G)$.  Let $g\in N_i$ be an element which is not fixed by $\phi$.  Then $\phi(g)=g\cdot z$, where $z\in Z_i$.  Since $\phi$ acts trivially on $G/\gamma_2^T(G)$ it acts trivially on $Z_i$ by Lemma \ref{l:trivial}.

Therefore, $\phi$ preserves the group generated by $g$ and $z$, which is abelian since the conjugation action of $g$ on $Z_i$ is trivial.  Therefore, $\langle g,z\rangle\cong \bZ^2$.  Choose an ordering on this copy of $\bZ^2$ which is not preserved by $\phi$.  By Rhemtulla's Extension Theorem, there exists an ordering on $N_i$ which restricts to the pre--chosen ordering on $\bZ^2$.  By Lemma \ref{l:soc}, we can pull this ordering back to $G$.  Since the ordering is not preserved on $N_i$, it is not preserved on $G$.
\end{proof}

In the remainder of this section we shall develop an alternative viewpoint on Theorem \ref{t:main} which makes the result more transparent, at least in the case of surface groups and free groups.  Recall that a finitely generated group $G$ is called {\bf hyperbolic}, {\bf Gromov hyperbolic} or {\bf negatively curved} if there is a $\delta\geq 0$ such that whenever $g,h\in G$, any geodesic in $G$ (with respect to the word metric) connecting $g$ and $h$ is contained in a $\delta$--neighborhood of the union of two geodesics connecting the identity to $g$ and $h$ respectively.  Being $\delta$--hyperbolic is a quasi--isometry invariant, though the precise value of $\delta$ which witnesses $\delta$--hyperbolicity depends on the generating set of $G$.

For basics on hyperbolic groups, the reader is referred to \cite{G}.  The property of hyperbolic groups we will be most interested in presently is the notion of the {\bf Gromov boundary} of an infinite hyperbolic group $G$, denoted $\partial G$.  Recall that to define $\partial G$, we fix a basepoint in $G$ and consider equivalence classes of geodesic rays emanating from the basepoint (in the Cayley graph of $G$).  Two geodesic rays are equivalent if they remain bounded distance from each other.  Using the $\delta$--hyperbolicity of $G$, it is possible to check that $\partial G$ is independent of the basepoint.

If two geodesic rays agree along long initial segments, then they are close.  It is easy to produce a dense set of points in $\partial G$ using the elements of $G$ itself.  Indeed, note that each infinite order element $g\in G$ gives rise to a point $x_g\in\partial G$ given by positive powers of $g$.  The precise statement is as follows, and a discussion can be found in \cite{G} (see also \cite{KapBenk}):

\begin{lemma}\label{l:hypback}
Let $G$ be a hyperbolic group.
\begin{enumerate}
\item
Each infinite order element $g\in G$ induces a loxodromic, fixed--point free isometry $\psi_g$ of the Cayley graph of $G$.
\item
For each infinite order $g$, the isometry $\psi_g$ has exactly two fixed points on $\partial G$, denoted $x_g$ and $y_g$.  These are the attracting and repelling fixed points of $\psi_g$ and are given by \[x_g=\lim_{n\to\infty}g^n\] and \[y_g=\lim_{n\to\infty}g^{-n}.\]
\item
The set of points $\{x_g\}$ for infinite order elements $g\in G$ is dense in $\partial G$.
\item
The maps $\bN\to G$ which sends $n$ to $g^n(b)$ is a quasi--isometric embedding for each basepoint $b$.
\item
If $g,h\in G$ do not generate an elementary (virtually cyclic) subgroup of $G$ then the fixed points of $g$ and $h$ on $\partial G$ do not coincide.
\item
If $G$ is torsion--free and $1\neq g\in G$, then there exists a unique $h\in G$ such that $g=h^m$ for some $m>0$ and $h$ is itself not a proper power.
\end{enumerate}
\end{lemma}

The points $\{x_g\}$ should be thought of as the {\bf rational points} in $\partial G$.  The motivation for this terminology is taken from lattices in $\bR^n$.  Notions akin to the Gromov boundary can be defined for non--negatively curved metric spaces, such as $\bR^n$.  From $\bR^n$ we obtain a natural boundary which is homeomorphic to $S^{n-1}$.  In this same way, the boundary of $\bZ^n$ should be thought of as $S^n$.  Then the rational points on the boundary are obviously given by lines through the origin, all of whose slopes are rational.

\begin{lemma}\label{l:boundary}
Let $G$ be a hyperbolic, residually torsion--free nilpotent, let $1\neq g\in G$ and let $\{\mP_{\al}\}$ be the set of positive cones on $G$ which contain $g$.  Let $h$ be the smallest root of $g$.  Then \[\bigcap_{\al}\mP_{\al}=\{h^n\mid n>0\}.\]
\end{lemma}
\begin{proof}
We must first check that this intersection is nonempty.  Clearly, $g$ is nontrivial in some torsion--free nilpotent quotient $N$ of $G$.  We may declare $g$ to be positive, thus defining a partial ordering on $N$.  By Rhemtulla's Extension Theorem we can extend this partial ordering to all of $N$, and then to all of $G$.  Therefore there is at least one positive cone which contains $g$.

Now suppose that $1\neq k\in G$ is another element such that $g$ and $k$ share no common powers.  Then there is minimal $i$ and a quotient $N_i=G/\gamma_i^T(G)$ in which $g$ and $k$ are both nontrivial.  Either the image of $\langle g,k\rangle$ in $N_i$ is isomorphic to $\bZ^2$ or it is cyclic.  If it is cyclic, replace $g$ and $k$ by powers so that they are equal in $N_i$.  We then take the smallest $j>i$ such that $g$ and $k$ are not equal in $N_j$.  It follows that $g$ and $k$ differ in $N_j$ by an element of $Z_j$, so that the image of $\langle g,k\rangle$ in $N_j$ is isomorphic to $\bZ^2$.

Now choose an ordering on the copy of $\bZ^2$ we have produced in which the image of $g$ is positive and the image of $k$ is negative.  By Rhemtulla's Extension Theorem, this ordering extends to all of $N_i$ (or $N_j$).  By the standard ordering construction of Lemma \ref{l:soc}, we can pull back the resulting ordering to $G$ in which $g$ is positive and $k$ is negative.  Therefore, \[k\notin\bigcap_{\al}\mP_{\al}.\]  It follows that if \[k\in\bigcap_{\al}\mP_{\al}\] then $k$ and $g$ share a common power.  The existence of $h$ with the desired properties follows from Lemma \ref{l:hypback}.
\end{proof}

It follows that the space $LO(G)$ recovers the Gromov boundary of a residually torsion--free nilpotent hyperbolic group in the following sense: the intersection of positive cones containing a given element of $G$ yields a unique point on the Gromov boundary of $G$, and the collection of all such points is a dense subset of $\partial G$.

Note that if $\phi\in\Aut(G)$ preserves each element of $LO(G)$ then it preserves the sets $\{g^n\mid n>0\}$ for elements $g$ which are not proper powers, since for each element of $G$, the automorphism $\phi$ preserves the positive cones which contain $g$.  

\begin{proof}[Proof of Theorem \ref{t:boundary}]
Suppose that $\phi\in\Aut(G)$ acts nontrivially on $LO(G)$.  Then there is a $g\in G$ and an ordering $P$ of $G$ such that $g\in P$ but $\phi(g)\notin P$.  It follows that $\phi(g)$ is not contained in the intersection of all positive cones in $G$ which contain $g$, so that $\phi(g)$ and $g$ share no common power.  It follows that $\phi(g)$ and $g$ cannot generate a virtually cyclic subgroup of $G$.

It follows that either $x_g$ and $x_{\phi(g)}$ are different in which case $\phi$ acts nontrivially on $\partial G$, or \[\lim_{n\to\infty}g^n=\lim_{n\to\infty}\phi(g^n).\]  If $x_{g^{-1}}$ and $x_{\phi(g^{-1})}$ are different then again we see that $\phi$ acts nontrivially on $\partial G$.  Therefore, we may assume that the quasi--geodesics determined by $g$ and $\phi(g)$ have the same endpoints at infinity.  But then the subgroup generated by $g$ and $\phi(g)$ is virtually cyclic, a contradiction.
\end{proof}
In connection with the proof of Theorem \ref{t:boundary}, we note the following:
Suppose $\phi\in\Aut(G)$ acts trivially on $LO(G)$.  Then $\phi$ preserves the sets of the form $\{g^n\mid n>0\}$ for elements $g$ which are not proper powers.  It follows that for each $1\neq g\in G$, the limits \[\lim_{n\to\infty}g^n\] and \[\lim_{n\to\infty}\phi(g^n)\] are equal, so that the rational points $\{x_g\}$ on $\partial G$ are preserved by $\phi$.  It follows that $\phi$ acts trivially on $\partial G$.

\section{Homology, orderings, residual finiteness and faithful representations}\label{s:homo}
In this section we will make some remarks about homology representations of $\Out(G)$, $O(G)$ and residual finiteness.  It would be nice if we could formulate and prove an analogous result to Theorem \ref{t:main} for the action of $\Out(G)$ on $O(G)$, but unfortunately we encounter various difficulties.  The proofs as they are given for $LO(G)$ will not work for $O(G)$.  One difficulty is the following: any residually finite group has a residually finite automorphism group.  On the other hand, it is not true that each residually finite group has a residually finite outer automorphism group.  In fact, Wise proves in \cite{W} that every finitely generated group embeds in the outer automorphism group of some residually finite group.

For certain residually torsion--free nilpotent groups however, it is possible to make $\Out(G)$ act faithfully on $O(G)$ just by exploiting the fact that the homology representation \[\Out(G)\to \Aut(H_1(G,\bQ))\] is faithful.  Consider $\Out(A_{\gam})$, where $\gam$ is a finite graph and $A_{\gam}$ is the associated right-angled Artin group.  Recall that $A_{\gam}$ is the free group on the vertices of $\gam$ together with the commutation relations between vertices whenever they are connected by an edge.  See Charney's expository article \cite{C} for more details.

Whereas abelian and free groups have very complicated automorphism groups, it is often the case that right-angled Artin groups have finite outer automorphism groups.  In fact, Charney and Farber have recently proved in \cite{ChFar} that a ``generic" right-angled Artin group has a finite outer automorphism group.

\begin{prop}
Let $A_{\gam}$ be a ``generic" right-angled Artin group.  Then $\Out(A_{\gam})$ acts faithfully on the abelianization $A_{\gam}^{ab}$ of $A_{\gam}$.  In particular, $\Out(A_{\gam})$ acts faithfully on $O(G)$.
\end{prop}
\begin{proof}
For a generic right-angled Artin group $A_{\gam}$, the outer automorphism group is generated by automorphisms of the graph $\gam$ and inversions of the vertices of $\gam$.  It follows easily that $\Out(A_{\gam})$ acts faithfully on the abelianization \[A_{\gam}^{ab}=A_{\gam}/\gamma_2^T(A_{\gam}).\]  For any given outer automorphism, one may choose an ordering on $A_{\gam}^{ab}$ which is not preserved by the action of that outer automorphism.  Any ordering on $A_{\gam}$ can be pulled back to an ordering on $A_{\gam}$ by the standard ordering construction of Lemma \ref{l:soc}.
\end{proof}

\section{Some final examples}
As claimed in the introduction, it is not true in general that $\Aut(G)$ acts faithfully on the left orderings $LO(G)$, nor is it true that $\Out(G)$ acts faithfully on the conjugacy classes in $LO(G)$ or on $O(G)$.  Consider, for instance, the fundamental group $K$ of the Klein bottle.  We have the presentation \[K=\langle x,y\mid x^{-1}yx=y^{-1}\rangle.\]  If $P$ is an ordering on $K$ then $P$ is certainly not bi--invariant.  Indeed, either $y\in P$ or $y\in P^{-1}$, but conjugation by $x$ takes $y$ to $y^{-1}$.

It is known that $K$ admits exactly four left--invariant orderings.  A discussion of this fact and other groups which admit only finitely many left--invariant orderings can be found in the book \cite{MR}.  It is easy to find various automorphisms of $K$ which have infinite order.  For instance, the automorphism which sends $x$ to $xy$ and fixes $y$ can easily be seen to have infinite order, whence it follows that $\Aut(K)$ is infinite.  Thus we see that $\Aut(K)$ cannot act faithfully on the space of left--invariant orderings $LO(K)$.

By the remarks above, we see that there are at most two conjugacy classes of left--invariant orderings on $K$, since if an ordering $P$ declares $y$ to be positive then a conjugate of $P$ declares $y$ to be negative.  It follows that if $\Out(K)$ acts faithfully on conjugacy classes of elements of $LO(K)$ then $\Out(K)$ can have at most two elements.  However:

\begin{prop}
$\Out(K)\cong\bZ/2\bZ\times\bZ/2\bZ$.
\end{prop}

One can check directly from the presentation of $K$ that the three non--inner automorphisms $\al_1:x\mapsto xy$, $\al_2:x\mapsto yx$, and $\al_3:x\mapsto x^{-1}$, where these are extended to $K$ by letting them fix $y$ in the first two cases and $\al_3:y\mapsto y^{-1}$, generate $\Out(K)$ and that they all have order two in $\Out(K)$.  Furthermore, $\al_1$ and $\al_2$ differ by an inner automorphism.

Another way to understand the outer automorphism group of $K$ is by an analogue of the Dehn--Nielsen--Baer Theorem, which shows that the mapping class group of the Klein bottle is actually $\bZ/2\bZ\times\bZ/2\bZ$.

The proof of Proposition \ref{p:klein} is now immediate.


\begin{thebibliography}{99}
\bibitem{A}J\o rgen Ellegaard Andersen.  Asymptotic faithfulness of the quantum ${\rm SU}(n)$ representations of the mapping class groups.  {\it Ann. of  Math}. 2, 163, no. 1, 347--368, 2006.
\bibitem{BL}Hyman Bass and Alexander Lubotzky.  Linear-central filtrations on groups.  {\it The mathematical legacy of Wilhelm Magnus: groups, geometry and special functions}.  Contemp. Math. 169, pp/ 45--98, 1994.
\bibitem{CB}Andrew J. Casson and Steven A. Bleiler.  {\it Automorphisms of surfaces after Nielsen and Thurston}.  London Mathematical Society Student Texts, 9, 1988.
\bibitem{C}Ruth Charney.  An introduction to right-angled Artin groups.  {\it Geom. Dedicata}, 125, 141--158, 2007.
\bibitem{ChFar}Ruth Charney and Michael Farber.  Random groups arising as graph products.  Preprint.
\bibitem{DDRW}Patrick Dehornoy, with Ivan Dynnikov, Dale Rolfsen and Bert Wiest.  {\it Ordering Braids}.  Mathematical Surveys and Monographs, volume 148.  American Mathematical Society, 2008.
\bibitem{G}Mikhail Gromov.  Hyperbolic groups.  {\it Essays in group theory}, 75--263, {Math. Sci. Res. Inst. Publ.}, 8, Springer, New York, 1987.
\bibitem{KapBenk}I. Kapovich and N. Benakli.  Boundaries of hyperbolic groups.  {\it Contemporary Mathematics}, vol. 296, 39--94, 2002.
\bibitem{KW}I. Kapovich and D. Wise.  Equivalence of some residual properties for word-hyperbolic groups.  {\it J. Algebra} 223, no. 2, 562--583, 2000.
\bibitem{K}Thomas Koberda.  Asymptotic linearity of the mapping class group and a homological version of the Nielsen-Thurston classification.  To appear in {\it Geom. Dedicata}.
\bibitem{M}A.I. Mal'cev.  On the full ordering of groups.  {\it Trudy Mat. Inst. Steklov.} 38, 173--175, 1951.
\bibitem{MR}Roberta Botto Mura and Akbar Rhemtulla.  {\it Orderable groups}.  Lecture notes in pure and applied mathematics, volume 27.  1977.
\bibitem{N}A. Navas.  On the dynamics of (left) orderable groups.  Preprint.  Available at arXiv: math.GR/0710.2466, 2007.
\bibitem{NR}A. Navas and C. Rivas.  Describing all bi--orderings on Thompson's group $F$.  {\it Groups Geom. Dyn.} 4, no. 1, 163--177, 2010.
\bibitem{R}M.S. Raghunathan.  {\it Discrete subgroups of Lie groups}.  Springer-Verlag, New York-Heidelberg, 1972.
\bibitem{Rh}Akbar Rhemtulla.  Right-ordered groups.  {\it Canad. J. Math.} 24, 891--895, 1972.
\bibitem{S}Adam S. Sikora.  Topology on the spaces of orderings of groups.  {\it Bull. London Math. Soc.} 36, no. 4, 519--526, 2004.
\bibitem{W}Daniel T. Wise.  A residually finite version of Rips' construction.  {\it Bull. London Math. Soc.} 35, 1, 23--29, 2003.
\end{thebibliography}
\end{document}